\documentclass[12pt]{article}
\usepackage{amsmath}
\DeclareMathOperator{\card}{card}
\begin{document}
\newtheorem{thm}{Theorem}[section]
\newtheorem{lem}[thm]{Lemma}
\def\N{{\mathbf N}}
\def\Q{{\mathbf Q}}
\def\R{{\mathbf R}}
\def\C{{\mathbf C}}
\def\P{{\mathbf P}}
\def\Z{{\mathbf Z}}
\def\v{{\mathbf v}}
\def\O{{\mathcal O}}
\def\M{{\mathcal M}}
\def\tr{\mbox{Tr}}
\def\qed{{\tiny $\clubsuit$ \normalsize}}

\Large {\bf Counting Rational Points on Ruled Varieties}\normalsize

\vspace{.2in}

\begin{center}
\begin{tabular}{l}
David McKinnon \\
Department of Pure Mathematics \\
University of Waterloo \\
Waterloo, ON\ \  N2L 3G1 \\
Canada \\
Email: dmckinnon@math.uwaterloo.ca \\
\end{tabular}
\end{center}

\begin{abstract}
In this paper, we prove a general result computing the number of
rational points of bounded height on a projective variety $V$ which is
covered by lines.  The main technical result used to achieve this is
an upper bound on the number of rational points of bounded height on a
line.  This upper bound is such that it can be easily controlled as
the line varies, and hence is used to sum the counting functions of
the lines which cover the original variety $V$.
\end{abstract}


\section{Introduction}

In algebraic geometry, the general notion of studying an algebraic variety
by studying families of curves which cover it is a very old and fruitful
one.  However, it has not been much used to study the density of rational
points on algebraic varieties, because there has not been the necessary
uniformity in the results for lower-dimensional varieties and their 
counting functions.

Heath-Brown suggests in [H-B] that the technique could be quite widely
used, and gives an example of how to compute the counting function of
a certain cubic surface by studying families of cubic curves lying on
it.  In the spirit of this idea, we will describe in this paper a
method for obtaining upper bounds (and in many cases, asymptotic
formulae) for the counting functions for rational points with respect
to height on projective varieties which are covered by lines.

In particular, we consider an algebraic variety $V$ defined over a
number field $K$, embedded in projective space $\P^n$ for some $n$,
and consider the usual (multiplicative) height function, not
normalised to be independent of the field $K$.  We define the counting
function thus:
\[N_L(B) = \card\{P\in V(k) \mid H(P)\leq B\}\]
The counting function counts the number of $K$-rational points of $V$
whose height is at most $B$.  (That this is well defined is an
immediate consequence of a theorem of Northcott -- see for example
[Vo], Proposition 1.2.9.(g).)

If $V$ is covered by a set $\mathcal{L}$ of lines, then we might hope
to compute the counting function of $V$ by summing the respective 
counting functions for the lines in $\mathcal{L}$:
\[N_V(B) = \sum_{L\in\mathcal{L}} N_L(B)\]
There are several excellent estimates for $N_L(B)$ in the literature.
The first of these was that of Schanuel [Sch], in which Schanuel
calculates very precisely the asymptotics of the counting function for
$\P^1$ over a number field $K$.  The chief drawback of this is that
the constant in the leading term depends on the specific embedding of
the line into $\P^n$, so for our purposes we will need a more specific
calculation.

This, too, has been done, by Thunder [Th], in which he calculates the
asymptotics of the counting function for an arbitrary line in $\P^n$
over an arbitrary number field.  Thunder makes clear that the leading
term in the counting function for a line $L$ is:
\[\frac{c_K}{H(L)} B^2\]
where $H(L)$ denotes the height of the Pl\"ucker point corresponding
to $L$, and $c_K$ is a constant depending only on the number field
$K$.  However, for our purposes, we will want strict control of how
many points of small height lie on lines of large height, so the
presence of an unbounded error term of any kind (which is definitely
necessary in all theorems of the sort that Schanuel and Thunder
sought) is fatal to our line of reasoning.

Thus, we must prove yet another result about heights of rational points
on lines, one which allows such a strict control -- this result is
Theorem~$\ref{main}$.  To get this control, we sacrifice the quality of
the constant $c_K$, and we obtain only an upper bound, rather than a 
lower one.  However, Theorem~$\ref{main}$ is still good enough to
give exactly the right exponent on $B$ in the counting function for
many algebraic varieties (see Theorem~$\ref{bundle}$), so our sacrifices
are certainly outweighed by our gains.

Theorem~$\ref{bundle}$ fits into the extensive literature which computes
the counting functions of algebraic varieties, which is too large to
summarise in a satisfactory fashion here.  Suffice it to say that it
is compatible with the conjectures of Batyrev and Manin, and that a 
more comprehensive overview of the history of the subject can be found
in [Si].

\section{Rational Points on Lines}

Let $K$ be an algebraic number field with ring of integers $\O_K$, and let
$L$ be a line in $\P^n_K$.  We wish to compute an upper bound for the
counting function:
\[N_L(B) = \card\{P\in L \mid H(P)\leq B\}\]
where $H(P)$ denotes the standard height in projective space:
\[H([x_0:\ldots :x_n]) = \prod_{v} \max_i\{|x_i|_v\})\] 
where $v$ ranges over all (isomorphism classes of) valuations on $K$.
Note that we do not normalise the height to be independent of the
field $K$.

Schanuel [Sch] derived some very precise asymptotics for $N_L(B)$:
\[N_L(B) = cB^2 + E(B)\]
where $c$ is a specific constant and $E(B)$ is an error term which is
$o(B^2)$.  Schanuel computes both of these quite precisely in the case
that $L=\P^1_K$.

However, this estimate will not suffice for our present purpose, since
we wish to control the set of points of small height on our lines as
well.  Furthermore, since our lines will not generally be identical to
$\P^1$, we wish to explore the dependence of $c$ on the height $H(L)$
of the line $L$, which we define to be the height of the corresponding
Pl\"ucker point in the Grassmannian $G(1,n)$.  An asymptotic version of
this has been derived by Thunder [Th], but like Schanuel's result, does
not control the behaviour of the points of small height.  

Thus, say $L$ corresponds to a 2-dimensional subspace of $K^{n+1}$,
spanned by the vectors $(a_0,\ldots,a_n)$ and $(b_0,\ldots,b_n)$.
We define:
\[H(L) = H((a_0dx_0+\ldots +a_ndx_n)\wedge (b_0dx_0 + \ldots + b_ndx_n))\]
where the result of the wedge product is interpreted as a point in
$\P^{(n^2+n)/2}_K$ with homogeneous coordinates $\{dx_i\wedge dx_j\}$
for $i\neq j$.

We can now state the main result of this section:

\begin{thm}\label{main}
The counting function for $L$ satisfies the following inequalities:
\[N_L(B) \leq \frac{c_K}{H(L)}B^2 + 1\]
where $c_K$ is a positive real constant depending only on the field $K$.
\end{thm}

\noindent
{\it Proof:} \/ Our first step will be to identify the $K$-rational
points of $L$ with a set of lattice points in a finite-dimensional
Euclidean space.  Let $[x_0\colon \ldots\colon x_n]$ be a $K$-rational
point on $L$.  By clearing denominators, we can ensure that
$x_i\in\O_K$ for all $i$.  In fact, by choosing a fixed set
$\mathcal{J}$ of representatives for the class group of $K$, we can
ensure that the coordinates $x_i$ generate an ideal in $\mathcal{J}$.  
This representation for $[x_0\colon\ldots\colon x_n]$ is unique up to
multiplication by a unit of $\O_K$.

Let $M$ be the rank two $\O_K$-module $M$ in $K^{n+1}$ consisting of
all the vectors in $L$ whose coordinates all lie in $\O_K$.  Let
$d=[K:\Q]$, and denote by $\sigma_1,\ldots ,\sigma_{r_1}$ the embeddings
of $K$ into $\R$, and by $\tau_1,\ldots ,\tau_{r_2}$ the embeddings of 
$K$ into $\C$, where $d=r_1+2r_2$.  

Using these embeddings, we can embed $M$ as a lattice of rank $2d$ in
$V=(\R^{n+1})^{r_1}\oplus(\C^{n+1})^{r_2}$.  We will abuse notation by
hereafter identifying $M$ with its image in $V$.  Define:
\[|(a_0,\ldots,a_n)|_i = \left\{\begin{array}{ll}
\max_j(|\sigma_i(a_j)|) & \mbox{if $i\leq r_1$} \\
\max_j(|\tau_{i-r_1}(a_j)|^2) & \mbox{if $i>r_1$}
\end{array}\right.\]
Thus, we can view $[x_0\colon\ldots\colon x_n]$ as a point in $M$
whose coordinates generate an element of the fixed set $\mathcal{J}$ of
ideals.  Thus, there exists a positive constant $c_1$ depending only on
the field $K$ such that:
\begin{equation}\label{height}
\frac{1}{c_1}\prod_i |(x_0,\ldots,x_n)|_i \leq H([x_0\colon\ldots\colon
x_n]) \leq c_1\prod_i |(x_0,\ldots,x_n)|_i
\end{equation}
We also define:
\begin{equation}\label{norm}
\|(x_0,\ldots,x_n)\| = \max_i\{|(x_0,\ldots,x_n)|^{d_i}_i\}
\end{equation}
where $d_i$ is 1 or $1/2$, depending on whether $i$ corresponds to a 
real or complex embedding, respectively.

If $\epsilon\in\O_K^*$ is a unit, then we have the relation:
\[|\epsilon(x_0,\ldots,x_n)|_i = |\epsilon|_i|(x_0,\ldots,x_n)|_i\]
where $|\epsilon|_i$ represents the absolute value of $\epsilon$ with
respect to the embedding $\sigma_i$ (if $i$ is at most $r_1$) or
$\tau^2_{i-r_1}$ (if $i$ is greater than $r_1$).  Thus, by the
Dirichlet Unit Theorem (see for example [Ne], Theorem 7.3), there is a
positive real constant $c_2$ depending only on $K$ such that for any
element $\v=(v_1,\ldots,v_{r_1+r_2})\in \R^{r_1}\oplus\C^{r_2}$, there
exists a unit $\epsilon\in\O_K^*$ such that $|\epsilon v_i|^{d_i}\leq
c_2|\epsilon v_j|^{d_j}$ for all $i$ and $j$, where $d_i$ and $d_j$ are
as in $(\ref{norm})$.

Applying this to $\v =
(|(x_0,\ldots,x_n)|_1,\ldots,|(x_0,\ldots,x_n)|_{r_1+r_2})$ reveals
that there is a positive real constant $c_3$ depending only on $K$
such that through multiplication by a suitable unit, we may assume
that for all $i$ and $j$:
\[|(x_0,\ldots,x_n)|^{d_i}_i\leq c_3|(x_0,\ldots,x_n)|^{d_j}_j\]

Thus, by equations $(\ref{height})$ and $(\ref{norm})$, we may find a
positive real constant $c_4$ depending only on $K$ such that for all
$K$-rational points $P\in L$, we can find a representation
$P=[x_0\colon\ldots\colon x_n]$ as above such that:
\[\frac{1}{c_4}H(P)\leq\|(x_0,\ldots,x_n)\|^d\leq c_KH(P)\]
Thus, when calculating an upper bound for $N_L(B)$, it will suffice to
compute an upper bound for the following function:
\[N'_L(B) = \{\v=(x_0,\ldots,x_n)\in M'\mid\|\v\|\leq B\}\]
where $M'$ denotes the set of vectors in $M$, counted modulo the
action of $K^*$.  In particular, we have:
\[N'_L(B) \geq N_L(B)^d\]

We will use the following well known result (it follows, for example,
from work of Thunder [Th]):

\begin{lem}
There is a positive real constant $\alpha$ depending only on $K$ such
that the determinant of $M$ is equal to $\alpha H(L)$.
\end{lem}

We therefore have reduced to showing that there is a positive real
constant $c$, depending only on the field $K$, such that:
\[N'_L(B)\leq \frac{c}{\det(M)}B^{2d} + 1\]

Thus, fix a positive real number $B$.  If $N_L(B)\leq 1$, then the
result is clear, so assume that $N_L(B)\geq 2$.  Then we can find two
$K$-linearly independent lattice points $P_1$ and $P_2$ in $M$ with
$H(P_i)\leq B$, since $K$-linearly dependent points in $M$ contribute
only one point to $N_L(B)$.

Choose a basis $\{\alpha_1,\ldots,\alpha_d\}$ for $\O_K$ over $\Z$,
and consider the set:
\[\{\alpha_1P_1,\ldots,\alpha_dP_1,\alpha_1P_2,\ldots, \alpha_dP_2\}\]
There is a positive real constant $\alpha$ depending only on $K$ such
that the height of $\alpha_j P_i$ is at most $\alpha H(P_i)$.

For any real number $H$, consider the set 
\[V_M(H)=\{\v\in M\otimes\R \mid \|\v\|\leq H\}\]
This set is convex and centrally symmetric, so by the previous
arguments, it follows that the real simplex spanned by the vectors
$\alpha_j P_i$ is contained in $V_M(\alpha B)$.  In particular, we
conclude that every element of $L$ with height at most $B$ corresponds
to a point of $M$ which is a vertex of a $2d$-dimensional real simplex
which is entirely contained in $V_M(\alpha B)$ and whose vertices are
all elements of $M$.

The number of such simplexes is at most
$\frac{V\alpha^{2d}}{\det(M)}B^{2d}$, where $V$ is the volume of the
standard $2d$-dimensional real simplex.  Each such simplex has $2d+1$
vertices, so we conclude that:
\[N'_L(B)\leq \frac{V(2d+1)\alpha^{2d}}{\det(M)}B^{2d}+1\]
and hence the theorem follows.  \qed

\section{Ruled Varieties}

In the spirit of Heath-Brown's remark in [H-B], Theorem $\ref{main}$
enables us to give easy upper bounds for the counting functions for
rational points on ruled varieties.  For instance, consider the following
situation.

Let $V\subset\P^n$ be a projective variety defined over a number field
$K$.  Assume that $V$ admits a $K$-rational morphism $\phi\colon
V\rightarrow X$ to a projective variety $X$ over $K$ such that the
fibres of $\phi$ are lines.  Then we can define a morphism $\psi\colon
X \rightarrow G(1,n)$ by $\psi(P) = [\phi^{-1}(P)]$, where $G(1,n)$
denotes the Grassmannian of lines in $\P^n$.  Let $D$ be the Pl\"ucker
divisor on $G(1,n)$ -- that is, the pullback of $\O(1)$ via the
Pl\"ucker embedding of $G(1,n)$.  Note that $\psi$ is injective, since
$\phi$ is a morphism.  We now have the following result:

\begin{thm}\label{bundle}
Using the notation of the previous paragraph, assume that the counting
function of $X$ with respect to the divisor $A=\psi^*(D)$ satisfies:
\[N_X(B) = \card\{P\in X(K) \mid H_A(P)\leq B\} = O(B^\epsilon)\]
for some $\epsilon<1$.  Then we have:
\[(1/c)B^2\leq N_V(B)\leq cB^2\]
for some positive constant $c$.
\end{thm}

\noindent
{\it Proof:} \/ The first inequality is clear by Schanuel's Theorem
[Sch], since $V$ contains at least one $K$-rational line.  Thus, we
turn our attention to the second inequality.  Write $H$ for the usual
height function in $\P^n$, and let $F=\phi^*A$.  Via the height
machine (see for example [Vo], Proposition 1.2.9), we obtain a
constant $\alpha$ such that for all $K$-rational points $P$ of $V$:
\begin{equation}\label{hart}
H_F(P)\leq \alpha H(P)
\end{equation}

We can now calculate as follows:
\begin{eqnarray*}
N_V(B) & \leq & \sum_{P\in X(K),\, H_A(P)\leq B} \alpha N_{\phi^{-1}(P)}(B) \\
& \leq & \sum_{P\in X(K),\, H_A(P)\leq B} \alpha\left(\frac{c_K}{H_A(P)}B^2 
+ 1\right)
\end{eqnarray*}
where this last inequality is by Theorem~$\ref{main}$ and the fact
that $H(\phi^{-1}(P)) = H_A(P)$.  The hypothesis of the theorem now
easily implies that this sum is asymptotically less than $cB^2$ for a
positive constant $c$, and the theorem is proven.  \qed

\vspace{.1in}

\noindent
{\bf Remarks:} In particular, Theorem~$\ref{bundle}$ applies to all
(relatively) minimal ruled surfaces (see section V.2 of [Ha] for a
discussion of such surfaces).  (This is not quite true, since the two
cases of $\P^1\times\P^1$ and $\P^2$ blown up at a single point do not
satisfy the hypotheses of Theorem~$\ref{bundle}$, but they can be 
handled in a similar manner, or indeed by any number of elementary 
approaches as well.)

The arithmetic of relatively minimal ruled surfaces over a rational
base curve has been dealt with admirably in several places, including
most notably in the very general treatment of Batyrev and Tschinkel
[BT] in the context of toric varieties, and in a more specific and
explicit way by Billard [Bi].  Note that in both these works, the
authors not only obtain the exponent in the leading term of the
counting function, but they also compute the constant in the leading
term and compute error terms, neither of which we are able to do here.

Finally, we remark that the results of Theorem~$\ref{bundle}$ are
consistent with the conjectures of Batyrev and Manin [BM].

\vspace{.1in}

Theorem~$\ref{main}$ can in principle be applied to any variety which 
is a union of lines in $\P^n$, by the simple expedient of summing the
counting functions of the individual lines, and controlling the point
of smallest height on each line.  Such an analysis proceeds trivially
for $\P^n$, for example, which is the union of a pencil of lines through
a fixed point, and the exponent on the upper bound thereby obtained is
sharp ($n+1$).  Similar analyses can be done for cones -- in both cases,
the point of smallest height on (almost all) lines is the basepoint of
the linear system which sweeps out the variety.

One might hope to obtain analogues of Theorem~$\ref{main}$ for curves
other than lines.  Indeed, Heath-Brown in [H-B] did so with his
Theorems 2 and 3 (the latter is conditional on a certain hypothesis on
the ranks of elliptic curves).  However, the chief advantage of 
Theorem~$\ref{main}$ is that the asymptotic growth of the counting 
function shrinks as the height of the line grows, making it much
easier to sum counting functions over infinitely many lines.
It would be interesting to try to obtain analogues of Theorem~$\ref{main}$
for curves of higher degree.

\end{document}